\documentclass{article}
\usepackage{amssymb,latexsym,amsmath,amsthm,mathrsfs}
\usepackage{graphicx}
\newcommand{\comment}[1]{}

\begin{document}
\title{An inquiry into whether or not 1000009 is a prime number\footnote{Presented to the St. Petersburg Academy on March 16, 1778.
Originally published as
{\em Utrum hic numerus 1000009 sit primus necne inquiritur},
Nova acta academiae scientiarum Petropolitanae \textbf{10} (1797),
63--73.
E699 in the Enestr{\"o}m index.
Translated from the Latin by Jordan Bell,
Department of Mathematics, University of Toronto, Toronto, Ontario, Canada.
Email: jordan.bell@gmail.com}}
\author{Leonhard Euler}
\date{}
\maketitle

\maketitle

1. Since this number is manifestly a sum of two squares,
namely $1000^2+3^2$, the question becomes this:
can this number be divided into two squares in any other way?
For if this cannot be done in any way, the number will certainly be
prime;\footnote{Translator: In Prop. 6, \S 35 of his 1749 {\em De numeris qui sunt aggregata duorum quadratorum}, E228, Euler shows that if $N=a^2+b^2$
where $a$ and $b$ are relatively prime and $N$ has no more representations
as a sum of two squares, then $N$ is prime.}
on the other hand, if this resolution could be done in some other way then
it will certainly not be a prime number, and its divisors could even be assigned.
Thus if we set one of the squares $=xx$, it needs to be inquired
whether the other one, namely $1000009-xx$, can escape being a square, except for the cases $x=3$
and $x=1000$. This can be investigated in the following way.\footnote{Translator:
The idea of this paper is that if some number is a sum of two squares
in two ways, then some other smaller number must be a square. We then check
all the cases. If we find a case where this smaller number is a square then
we can use this to find a factor of the original number, which is therefore
composite. If we check all the cases and none of them are 
squares, then the original number can be written as a sum of two squares in
only one way and hence is prime.}

2. Since the number ends in $9$, one of the squares is necessarily divisible
by $5$, and indeed thus by $25$.\footnote{Translator: Since
$1000009 \equiv 9 \pmod{10}$,
either $1000009-xx,xx \equiv 0,9 \pmod{10}$ or $1000009-xx,xx \equiv 4,5 \pmod{10}$.}
Let us therefore take the formula $1000009-xx$ to be divisible by $25$,
and it is clear that it necessarily happens that $x=25a+3$;
then this formula will be obtained:
\[
1000000-6\cdot 25a-25^2aa,
\]
which divided by $25$ becomes $40000-6a-25aa$. This form therefore should be
a square.\footnote{Translator: Since $1000009-xx$ is a multiple
of $25$ and they are both squares, $(1000009-xx)/25$ is a square.}

3. At this point two cases must be considered,
according to whether
$a$ is an even or an odd number. In the first case let $a=2b$, and
by dividing by 4, this resulting formula must also be a square:
\[ A=10000-3b-25bb. \]
For the other case let $a=4c+1$, and the resulting form will be the square
\[ B=39969-224c-400cc, \]
which is at any rate able to be an odd square;\footnote{Translator: I would guess that Euler
just means that there is no obvious reason why $B$ can't be a square,
not that it obviously can be a square.} on the other hand
for the same case let us take $c=4d-1$, and this formula results:
\[ C=39981+176d-400dd, \]
which, since it leaves $5$ when divided by $8$,
can never be a square.\footnote{Translator: Since the squares modulo $8$
are $0,1,4$.}
Therefore 
only the two formulas $A$ and $B$ need to be examined.

\begin{center}
{\Large The expansion of the formula
\[B=39969-224c-400cc.\]}
\end{center}

4. Here let us successively take all the values $0,1,2,3$ etc., both
positive and negative, for the letter $c$. Since 
the formula $400cc \pm 224c$ needs to be subtracted from the absolute
number $39969$, according to $c$ being a positive or negative
number, 
let us record these numbers successively subtracted in two columns,
along with their differences:

\begin{tabular}{p{1cm}|p{3cm}|p{1cm}|p{1cm}|p{3cm}|p{1cm}}
$c$&$400cc-224c$&Diff.&$c$&$400cc+224c$&Diff.\\
\hline
0&0&&0&0&\\
&&176&&&624\\
1&176&&1&624&\\
&&976&&&1424\\
2&1152&&2&2048&\\
&&1776&&&2224\\
3&2928&&3&4272&\\
&&2576&&&3024\\
4&5504&&4&7296&
\end{tabular}

It's clear right away here that in both cases the differences
continually increase by $800$.

5. These differences are then continuously subtracted from the absolute
number 39969; for convenience this is done in two columns, so that
it can be seen whether the numbers that result from this are squares:

\begin{minipage}{\textwidth}
\begin{tabular}{p{2cm}|p{2cm}||p{2cm}|p{2cm}}
39969&39969&31089&28849\\
176&624&4176&4624\\
\hline
39793&39345&26913&24225\\
976&1424&4976&5424\\
\hline
38817&37921&21937&18801\\
1776&2224&5776&6224\\
\hline
37041&35697&16161&12577\\
2576&3024&6576&7024\\
\hline
34465&32673&9585&5553\\
3376&3824&7376&\\
\hline
31089&28849&$* 2209$&\\
\end{tabular}
\end{minipage}

6. In both sides of the calculation, the single square that occurs is
$* 2209=47^2$;
whence it it apparent that the given number is not prime;
although this number is included in the paper
{\em De tabula numerorum primorum
usque ad millionem et ultra continuanda}, included in volume XIX of our
Novi Commentarii, 
it has divisors; for finding these it can be noted that
the square arose from the value $c=-10$, whence $a=-39$;
then collecting, $x=25a+3=-972$, and hence
\[ 1000009-xx=55225=235^2, \]
so that we thus have this double resolution:
\[ 1000^2+3^2=972^2+235^2, \]
and then by transposing
\[ 1000^2-235^2=972^2-3^2, \]
from which it follows
\[ (1000-235)(1000+235)=(972-3)(972+3), \]
or $1235\cdot 765=969\cdot 975$.
Then it would be $\frac{1235}{975}=\frac{969}{765}$, and these fractions are reduced
to this simplest one, $\frac{19}{15}$, from which it is then
concluded\footnote{Translator: See Prop. 7 of E228.
Say $N=a^2+b^2=c^2+d^2$, with both $b$ and $d$ odd.
Then $a^2-c^2=d^2-b^2$ and so $(a-c)(a+c)=(d-b)(d+b)$.
Let $k=\gcd(a-c,d-b)$. Then for some $l,m$ with $\gcd(l,m)=1$,
$a-c=kl$ and $d-b=km$, and so
$l(a+c)=m(d+b)$. Because $l$ and $m$ are relatively prime, 
$m$ divides $a+c$, so $a+c=mn$ for some $n$. Hence $ln=d+b$. Then it turns
out that $N=(\frac{k^2+n^2}{2})(\frac{m^2+l^2}{2})$;
we can check this by expanding this, and we end up getting
$\frac{a^2+b^2+c^2+d^2}{2}$, which indeed is $=N$. In the case
Euler is doing here, we get $l=15,m=19,k=51,n=65$.}
that our number has a common divisor with the sum of squares
$19^2+15^2$, which will thus be $293$, and we find that
\[
1000009=293\cdot 3413.
\]
From this it is clear that an error crept into the table in the above
mentioned paper, where all the prime numbers contained between
$1000000$ and $1002000$ were listed, which happened as
examination of the prime divisor $293$ was overlooked.

\begin{center}
{\Large The expansion of the formula
\[A=10000-3b-25bb\].}
\end{center}

7. This formula is a hundredth part of the formula $1000009-xx$, and for
its expansion
again two cases need to be distinguished, one in which $b$ is an even
number and another in which it is odd. It is evident that in the first
case, unless $b$ is an evenly even number the proposed formula cannot be a square.
Therefore let $b=4c$, and the resulting form, divided by $4$, will be
$2500-3c-100cc$. It can be seen without too much difficulty that this cannot
be a square aside from the case $c=0$.\footnote{Translator: The case $c=0$ just gives $1000009=1000^2+3^2$, and we are looking for other decompositions into
two squares.}
First it's evident that it can't happen when $c= \pm 1$; next, it likewise cannot
happen either when $c= \pm 2$.\footnote{Translator: I don't see why it's
obvious that $c=\pm 1$ doesn't give squares. But there the numbers are $2397$
and $2603$, and we can just check that these are different from the squares
$48^2,49^2,50^2,51^2,52^2$. For $c=\pm 2$ they are not squares because neither
are divisible by $4$.}
Thus let $c= \pm 3$, and our formula turns into $2500-900 \pm 9=1600 \pm 9$,
which cannot be a square. Now if one takes $c= \pm 4$, it is
\[
2500-1600 \pm 12 = 900 \pm 12,
\]
which is certainly not a square. Next, even taking $c= \pm 5$ a square
still doesn't arise; for this yields
\[
2500-2500 \pm 15=0 \pm 15.
\] 

8. For the other case where $b$ is an odd number,
one first puts $b=4d+1$,
and the proposed formula becomes
\[
9972-212d-400dd,
\]
which divided by $4$ will be
\[
2493-53d-100dd,
\]
which in the case $d=0$ is apparently not a square. 
Therefore let us take $d= \pm 1$, which gives $2393 \pm 53$,
also not a square;
and the case $d= \pm 2$ yields
$2093 \pm 106$;
indeed the case $d= \pm 3$ gives $1593 \pm 159$,
from neither of which result squares, nor for the case
$d= \pm 4$, which obviously gives $893 \pm 212$.
Next, the case $d=-5$ yields $-7+265$. Finally let
$b$ be a number of the form $4d-1$, which yields
\[
9978+188d-400dd.
\]
Since this number is even but not divisible by $4$,
it cannot be a square.

9. So that the strength of this calculation can be better shown,
let us examine the resolvability of another number of this kind into
two squares, which we shall take to be $1000081=1000^2+9^2$,
and let us see whether it can still be resolved into two squares
in another way. One of these, like in the preceding case,
is necessarily divisible by $5$.\footnote{Translator: The squares modulo $10$
are $0,1,4,5,6,9$. Either $1000081-xx,xx \equiv 5,1 \pmod{10}$,
$\equiv 5,6 \pmod{10}$ or $\equiv 5,9 \pmod{10}$.}
Therefore by putting one square $=xx$, let us see whether
the remaining part $1000081-xx$ can be a square divisible by 5, i.e. 25.

10. To this end let us set $x=25y+9$, and it will become this formula
\[
1000000-18\cdot 25y-25^2yy,
\]
which divided by 25\footnote{Translator: The original has the misprint 4.}
turns into this simpler one:
\[
40000-18y-25yy.
\]
Now first let $y$ be an even number, that is $y=2a$,
and dividing the formula again by $4$ one obtains:
\[
A=10000-9a-25aa.
\]
Second for the odd number let us put:
$1^\circ$ $y=4b+1$, and this gives
\[
B=39957-272b-400bb.
\]
This number is odd and leaves 5 when divided by 8, and hence cannot
be a square, whence the formula $B$ can be omitted altogether.
$2^\circ$ let us put $y=4c-1$, and the formula will be:
\[
C=39993+128c-400cc,
\]
where the number 39993 leaves 1 when divided by 8, and hence it is
appropriate to subject it to further investigation.

\begin{center}
{\Large The decomposition of the formula \[C=39993+128c-400cc\].}
\end{center}

11. Since here numbers contained in the form $400cc \pm 128c$
need to be successively subtracted from the absolute number $39993$,
as above to assist this calculation
we
write in the following table
the numbers to be subtracted along with their differences,
according to whether $c$ is positive or negative:

\begin{tabular}{|p{1cm}|p{3cm}||p{1cm}|p{1cm}|p{3cm}|p{1cm}}
$c$&$400cc-128c$&Diff.&$c$&$400cc+128c$&Diff.\\
\hline
0&0&&0&0&\\
&&272&&&528\\
1&272&&1&528&\\
&&1072&&&1328\\
2&1344&&2&1856&\\
&&1872&&&2128\\
3&3216&&3&3984&
\end{tabular}

Here again the differences successively increase by 800.

12. Thus let us subtract these differences continually increasing by eight hundred
from the absolute number 39993. The calculation goes as follows: 

\begin{tabular}{p{2cm}|p{2cm}||p{2cm}|p{2cm}}
39993&39993&30633&29353\\
272&528&4272&4528\\
\hline
39721&39465&26361&24825\\
1072&1328&5072&5328\\
\hline
38649&38137&21289&19497\\
1872&2128&5872&6128\\
\hline
36777&36009&15417&13369\\
2672&2928&6672&6928\\
\hline
34105&33081&8745&6441\\
3472&3728&7472&\\
\hline
30633&29353&1273&
\end{tabular}

Plainly no square occurs here.

\begin{center}
{\Large The decomposition of the formula \[A=10000-9a-25aa\].}
\end{center}

13. Let us put an even number in place of $a$,
which indeed must be evenly even,
and therefore let $a=4e$.
Therefore by dividing by 4 this form arises:
$2500-9e-100ee$. 
Then numbers contained in the form $100ee \pm 9e$ need to be successively
subtracted from the absolute number,
which are recorded in the following table,
according as $e$ is a positive or negative number:

\begin{tabular}{p{1cm}|p{3cm}|p{1cm}||p{3cm}|p{1cm}}
$e$&$100ee-9e$&Diff.&$100ee+9e$&Diff.\\
\hline
0&0&&0&\\
&&91&&109\\
1&91&&109&\\
&&291&&309\\
2&382&&418&\\
&&491&&509\\
3&873&&927&
\end{tabular} 

Thus let us continually subtract these differences increasing by
two hundred from the absolute number $2500$ in the following way:

\begin{tabular}{p{3cm}|p{3cm}}
2500&2500\\
91&109\\
\hline
2409&2391\\
291&309\\
\hline
2118&2082\\
491&509\\
\hline
1627&1573\\
691&709\\
\hline
936&864\\
891&\\
\hline
45&
\end{tabular}

where no squares occur besides 2500, which however leads to the known square
$1000^2$.

14. Now let $a$ be an odd number, and first of the form $4f+1$,
whence our formula turns into\footnote{Translator: The original
has the misprint $9966-236f-4ff$.}
\[
9966-236f-400ff.
\]
Since this number is oddly even, it cannot be a square.
Therefore let us put $a=4f-1$, and the formula becomes\footnote{Translator: The original has the misprint $9984+164-4ff$.}
\[
9984+164f-400ff,
\]
which is evenly even; and divided by 4 it turns into this:
\[
2496+41f-100ff.
\]
Therefore here numbers of the form $100 \pm 41f$ need to be subtracted
from the absolute number, which, according as $f$ is a positive or negative
number, works out like this:

\begin{tabular}{p{1cm}|p{3cm}|p{1cm}||p{1cm}|p{3cm}|p{1cm}}
$f$&$100ff-41f$&Diff.&$f$&$100ff+41f$&Diff.\\
\hline
0&0&&0&0&\\
&&59&&&141\\
1&59&&1&141&\\
&&259&&&341\\
2&318&&2&482&\\
&&459&&&541\\
3&777&&3&1023&
\end{tabular}

Now let these differences increasing by two hundred be subtracted
from the absolute number:

\begin{tabular}{p{3cm}|p{3cm}}
2496&2496\\
59&141\\
\hline
2437&2355\\
259&341\\
\hline
2178&2014\\
459&541\\
\hline
1719&1473\\
659&741\\
\hline
1060&732\\
859&\\
\hline
201&
\end{tabular}

Therefore because no square occurs anywhere in this calculation, it is certain
that the proposed number 1000081 can be resolved into two squares
in only one way, and hence it is certain that it is prime, as was
presented in the table of the above mentioned paper; and it is
particularly noteworthy that we have determined this  truth by such 
an easy calculation.

15. It is regrettable however that this method cannot be applied to exploring
all numbers, but is restricted just to those numbers which not only
are sums of two squares, but also which end in 1 or 9;
for so much of the success of this method is because one of
the squares is divisible by 5.

16. In the meantime however, plainly all numbers contained
in the form $4n+1$ ending in either 1 or 9 can be examined by this
method with equal success;
for we have found that if such a number can be resolved
into two squares, one of them is certainly
divisible by 5.
Then also, by calculating according to the rule that has been established,
if it turns out that the proposed number can be resolved into two squares
in just one way, this is certain proof that it is prime;
but on the other hand, if it can be done in two or more ways, from this
factors can be assigned as we did above.
Indeed if it happens that the proposed number cannot be divided
into two squares even in one way, then this also is proof that the
number is not prime, even if its factors cannot then be defined;
for it can be concluded that it has at least two prime factors of the form
$4n-1$. 

\end{document}